\documentclass{amsproc}

\usepackage{amsfonts, amsthm, amsmath, amssymb, color, bbm, hyperref, euscript}

\RequirePackage{mathrsfs} \let\mathcal\mathscr

\newtheorem{theorem}{Theorem}
\newtheorem*{thm}{Theorem}
\newtheorem{lemma}[theorem]{Lemma}
\newtheorem{corollary}{Corollary}
\newtheorem*{rem}{Remark}
\newtheorem*{eg}{Example}

\theoremstyle{definition}

\theoremstyle{remark}

\numberwithin{equation}{section}

\newcommand{\delete}[1]{}

\renewcommand{\phi}{\varphi}

\renewcommand{\AA}{\mathbb{A}}

\renewcommand{\leq}{\leqslant}

\renewcommand{\geq}{\geqslant}

\renewcommand{\bar}{\overline}

\renewcommand{\b}{\mathbf{b}}
\renewcommand{\a}{\mathbf{a}}

\renewcommand{\rho}{\varrho}

\DeclareMathOperator{\Mod}{mod} 
\renewcommand{\bmod}[1]{\,(\Mod{#1})}

\begin{document}

% \title[short text for running head]{full title}
\title[Character sums of composite moduli and hybrid subconvexity]{Character sums of composite moduli and hybrid subconvexity} 

%    Only \author and \address are required; other information is
%    optional.  Remove any unused author tags.

%    author one information
% \author[short version for running head]{name for top of paper}
\author{Roman Holowinsky}
\address{Department of Mathematics\\ The Ohio State University\\100 Math Tower\\231 West 18th Avenue\\Columbus, OH 43210\\USA}
\email{holowinsky.1@osu.edu}

\author{Ritabrata Munshi}
\address{School of Mathematics\\Tata Institute of Fundamental Research\\1 Homi Bhabha Road\\Colaba\\Mumbai 400005\\India}
\email{rmunshi@math.tifr.res.in}

\author{Zhi Qi}
\address{Department of Mathematics\\ The Ohio State University\\100 Math Tower\\231 West 18th Avenue\\Columbus, OH 43210\\USA}
\email{qi.91@osu.edu}

%    author two information

\thanks{The first author completed this work through the support of the NSF grant DMS-1068043. The second author was partly supported by Swarna Jayanti
Fellowship, 2011-12, DST, Govt. of India.}

%\subjclass[2000]{Primary }
%    The 2010 edition of the Mathematics Subject Classification is
%    now available.  If you are citing a classification from the
%    new scheme, use the following input coding instead.
%\subjclass[2010]{Primary }

%\date{}

\dedicatory{On the occasion of James W. Cogdell's 60th birthday.}

\begin{abstract}
	Let $M=M_1 M_2 M_3$ be the product of three distinct primes and let $\chi=\chi_1 \chi_2 \chi_3$ be a Dirichlet character of modulus $M$ such that each $\chi_i$ is a primitive character modulo $M_i$ for $i=1,2,3$. In this paper, we provide a $\delta$-symbol method for obtaining non-trivial cancellation in smooth character sums of the form $\sum_{n=1}^\infty \chi(n) W(n/N)$, with $N$ roughly of size $\sqrt M$ and $W$ a smooth compactly supported weight function on $(0, \infty)$. As a corollary, we establish hybrid subconvexity bounds for the associated Dirichlet $L$-function. 
\end{abstract}

\maketitle

\section{Introduction and main results}
\label{intro}
Let $L(s,\pi)$ be the $L$-function associated with an irreducible cuspidal automorphic representation $\pi$ with unitary central character.  Analysis of $L(s,\pi)$ leads to information about the arithmetic or algebraic structure associated with $\pi$ and non-trivial estimates for $L(s,\pi)$ in terms of its analytic conductor $Q(s,\pi)$ (for values of $s$ in the critical strip $0\leqslant \Re(s) \leqslant 1$) often result in non-trivial applications.  One classical problem, the subconvexity problem, is to establish a bound of the form
$$
L(s,\pi) \ll Q(s,\pi)^{1/4-\delta}
$$
for some $\delta>0$ when $\Re(s)=1/2$.  In general, for $\Re(s)=1/2$, one has the convexity bound
$$
L(s,\pi) \ll_\varepsilon Q(s,\pi)^{1/4+\varepsilon}
$$
while the Riemann Hypothesis for $L(s,\pi)$ would imply the Lindel\"{o}f Hypothesis
$$
L(s,\pi) \ll_\varepsilon Q(s,\pi)^{\varepsilon}.
$$
Although the convexity bound is far from the expected Lindel\"{o}f bound, any power saving in the conductor is often sufficient for applications.  For example, subconvexity for $L\left(\tfrac{1}{2}+it, \textnormal{Sym}^2 f\right)$ and $L \left(\tfrac{1}{2}, \textnormal{Sym}^2 f \times \varphi\right) $, where $t$ is a fixed real number, $f$ is a varying holomorphic eigencuspform (with $\textnormal{Sym}^2 f$ its symmetric square) and $\phi$ is a fixed Hecke-Maass eigencuspform for the modular group $\textnormal{SL}_2(\mathbb{Z})$, implies the Mass Equidistribution Conjecture (a holomorphic analogue of the Quantum Unique Ergodicity Conjecture \cite{QUE}) for $\textnormal{SL}_2(\mathbb{Z})$.  The subconvexity problem has thus received much attention in various settings recently, however, a general method of proof for all $\pi$ does not yet exist.

In a collection of works by the authors, see for example \cite{Munshi:circle1}, \cite{Munshi:circle2}, \cite{Munshi:circle3}, \cite{Munshi:circle4}, \cite{Ho-Mu}, and \cite{Ho-Mu-Qi}, several methods have been developed to investigate the subconvexity problem  particularly in the case of Rankin-Selberg convolution $L$-functions where multiple parameters are varying.  Such methods have led to a variety of hybrid subconvexity results, most recently demonstrating that subconvexity bounds are more readily obtained for $L(\tfrac{1}{2}, \textnormal{Sym}^2 f \times \varphi)$ when both $f$ and $\varphi$ are varying. Indeed, in \cite{Ho-Mu-Qi} the authors establish the following result.

\begin{thm}
	Suppose  $k > \kappa \geqslant  2$ are integers, with $ k$ even, $P$ is a prime, $f$ is a Hecke cusp form of weight $k$ for $\mathrm {SL} (2,\mathbb{Z})$,
	and $g$ a newform of weight $2\kappa$ and level $P$.Then we have 
	\begin{equation*} %\label{eq: final bound}
	L\left( \tfrac 1 2, \textnormal{Sym}^2 f\otimes g \right) 
	\ll _{\varepsilon, \kappa}
	\begin{cases}
	k^{  {13}/ {29} } P^{ {25} /{29} } (kP)^{\varepsilon},  & \text { if } P^{ {13}/ {64} } < k \leqslant P^{ 4 /{13} },\\
	\left( P + k^{ {13} /7} P^{ 3 / 7} \right) (kP)^\varepsilon,  & \text{ if } P^{ 4 /{13} } < k < P^{ 3 / 8}.
	\end{cases}
	\end{equation*}
	This bound  beats the convexity bound $k P^{ 3 / 4} (kP)^{\varepsilon}$ when $P^{ {13} /{64}  + \delta} < k < P^{ 3 / 8 - \delta}$ for some $0 < \delta <  {11} / {128}$.
\end{thm}

%\textcolor{red}{STATE OUR MAIN THEOREM FROM PREVIOUS PAPER}

Although more parameters are contributing to the complexity and analytic conductor of the $L$-function in such hybrid subconvexity problems, these situations are amenable to a larger collection of analytic tools and methods.  For example, if $\pi=f_1 \times f_2$ with each $f_i$ a holomorphic newform of varying level $N_i$ and $(N_1,N_2)=1$, then one has several natural ``families'' and ``sub-families'' of $L$-functions to which $L(s,\pi)$ might be associated.  In order to prove subconvexity for $L(s,f_1 \times f_2)$, one might choose to first study a moment average over a basis of newforms of level $N_1$, of level $N_2$, or average over both $N_1$ and $N_2$.  If, instead, only one of the levels is varying, then we immediately lose that additional degree of freedom.  

Such hybrid subconvexity problems therefore raise a question regarding structure and which family/moment of $L$-functions one should consider.  In order to establish subconvexity in the case of $\pi=f_1 \times f_2$ above, it was seen in \cite{Ho-Mu} that one should average over the larger level family when studying a second moment while one should average over the smaller level family, as in \cite{HT}, when studying a first moment.  If one were to study the first moment over the larger level family, then one obtains exact evaluations of the moment average rather than subconvexity (see for example \cite{MR07}, \cite{FW09}, \cite{Nelson:stable}).  Ultimately, the subconvexity problem boils down to having a sufficient number of points of summation relative to the conductor and complexity of the $L$-function one is considering, without having too many points of summation.  

In an attempt to better understand the underlying structure of such hybrid subconvexity results, we turn to the classical example of GL(1) convolutions.  Of course, one has Burgess' well known result for Dirichlet $L$-functions of a primitive character $\chi$ of modulus $M$ (\cite[Theorem 3]{Burgess3}),
\begin{equation*}
L \left(\tfrac 1 2 + i t, \chi \right) \ll_{\varepsilon, t} M^{3/16 + \varepsilon}.
\end{equation*}
But this does not close the subject.  For example, recent work by Mili\'{c}evi\'{c} \cite{Milicevic} on powerful moduli, improves on Burgess' bound for Dirichlet characters with moduli a sufficiently large power of a prime.

In this paper, we present a method for obtaining subconvexity results when the modulus of the Dirichlet character is ``moderately'' composite. Specifically, when the modulus is a product of three distinct primes $M=M_1 M_2 M_3$. Our method is an adaptation of the one presented in \cite{Munshi:circle2}.  Since we are dealing only with Dirichlet characters, the method becomes more transparent.  However, our main result is weaker compared to the Burgess bound.  As such, this paper does not prove any new result and one should view this work mainly as pedagogical.  Our method easily generalizes in the case of ``highly'' composite moduli.  With the availability of more factors, one has more options to design a ``conductor lowering'' mechanism.  Nevertheless, we feel that in this case, the $q$-analogue of the van der Corput method (see Theorem 12.13 of \cite{IK}) is much stronger.

Recall that Burgess' bound  for Dirichlet $L$-functions relies on the estimation of the character sum (\cite[Theorem 1]{Burgess2}, \cite[Theorem 2]{Burgess3}),
$$ \sum_{n=N+1}^{N+H} \chi (n) \ll_{\varepsilon} H^{1/2} N^{3/16 + \varepsilon}. $$
When $\chi$ is of prime modulus $p$, the proof of this bound in \cite{Burgess1, Burgess3}  features the application of an important estimate of Weil for $\sum_{x\in \mathbbm {F}_p^\star} \chi (f(x))$. 

For our purpose we shall consider the {smooth} character sum
\begin{equation} \label{2eq: smooth character sum}
S_{\chi} (N) := \sum_{n =1}^\infty \chi (n) W \left(\frac n N \right),
\end{equation}
where $W $ is  a smooth weight function on $(0, \infty)$ supported in the interval $[1, 2]$ and satisfying $W^{(j)} (x) \ll_j 1$. We shall obtain the following result on this smooth character sum. Interestingly, our proof also depends on a certain bound due to Deligne and Fu for $
\sum_{\substack{(x_1, x_2)\in \mathbb{F}_p^{\star 2} }}\chi (x_1 x_2) 
e\left(f(x_1, x_2)/p\right)
$ which is rooted in algebraic geometry over a finite field like Weil's bound.

\begin{theorem}\label{thm: character sum}
	Let $M_1, M_2, M_3$ be three distinct primes and set $M:=M_1M_2M_3$.  Let $\chi_i$ be a primitive character modulo $M_i$ and set $\chi := \chi_1 \chi_2 \chi_3$. 
	For $$M_1 \ll N \ll M_1 \min\big\{M_2^{2/3}, M_3^{2}\big\}$$ we have
	\begin{equation} \label{SchiNbound}
	S_{\chi} (N)  \ll_\varepsilon \left(M_2^{1/2} M_3^{1/2}+ {M_1^{1/4} M_2^{1/2}  } {N^{1/4}} + M_3^{1/2}{N^{3/4}}\right)M^\varepsilon.
	\end{equation}
\end{theorem}

As a corollary, we get the following hybrid subconvexity result.

\begin{corollary}\label{cor: subconvexity}
	Let $M_i$, $M$, $\chi_i$ and $\chi$, for $i = 1, 2, 3$, be given as in {Theorem~{\ref{thm: character sum}}}.  Set $\theta_i := \log{M_i}/\log{M}$ so that $\theta_1 + \theta_2 + \theta_3 = 1$.  If
	\begin{equation} \label{eq: relations between theta, 3}
	2 \delta \leq  \theta_1 \leq \tfrac{1}{2} - \delta, \ 
	\theta_3 \leq \tfrac{1}{4} - \delta , \ 
	\tfrac{1}{2} + 3 \delta \leq \theta_1 + 2 \theta_3
	\end{equation}
	for some $\delta > 0$, then 
	\begin{equation*}
	L \left(\tfrac 1 2, \chi \right) \ll_{\varepsilon} M^{\frac 1 4 - \frac \delta 2 + \varepsilon}.
	\end{equation*}
\end{corollary}

\begin{eg}
	When $(\theta_1,\theta_2, \theta_3)=\big(\tfrac{5}{12}, \tfrac{5}{12}, \tfrac{1}{6}\big)$, one may choose $\delta=\tfrac{1}{12}$ so that in this case
	$$
	L\left(\tfrac{1}{2}, \chi\right) \ll M^{\frac{5}{24}+\varepsilon}.
	$$
\end{eg}

\begin{rem}
	Of course, the set of triples $(\theta_1, \theta_2, \theta_3)$ for which one obtains non-trivial estimates for $S_\chi(N)$ {\rm(}and therefore subconvexity bounds for the corresponding set of $L\left(\tfrac{1}{2}, \chi\right)${\rm)} can be extended upon permuting the subscripts $1, 2, 3$.
\end{rem}

We shall see that a moment average will not be necessary in establishing Theorem \ref{thm: character sum}. Instead, the appropriate number of points of summation will be introduced directly via a $\delta$-symbol method which we describe in the next section.  A similar method may be found in \cite{Munshi:circle2}.  Furthermore, such a technique with similar arguments would establish subconvexity in the case of $\chi=\chi_1 \chi_2$ when $\chi_1=M_1$ and $\chi_2=M_2^2$ and seemingly in higher rank cases when the conductor is of an appropriate form. However, we do not yet see a blanket general structure to classify all situations in which such a $\delta$-symbol method would establish subconvexity.

\section{Preliminaries}

\subsection{Dirichlet $L$-functions and character sums}\label{sec: L and sum}

For a positive integer $M \geq 2$ let $\chi$ be a primitive Dirichlet character of modulus $M$. The Dirichlet $L$-function for $\chi$ is given by
$$L(s, \chi) = \sum_{n=1}^{\infty} \frac {\chi (n)} {n^s},$$
where the series converges for $\Re(s) > 0$.
\delete{Our goal is obtaining the subconvexity bound of the central value of the Dirichlet $L$-function for $\chi$, that is,
	\begin{equation}\label{2eq: subconvexity}
	L\left(\tfrac 1 2 , \chi \right) \ll_{\delta} M^{1/4 - \delta},
	\end{equation}
	for a certain positive constant $\delta$.} \delete{It is known that such a subconvexity bound can be deduced from the following bound of the character sum 
	\begin{equation*} %\label{2eq: character sum}
	\sum_{n =1}^N \chi (n) \ll_{\sigma, \varepsilon} N^{1-\sigma} ,
	\end{equation*}
	for $N \ll M^{1/2 + \varepsilon}$, where $ \sigma, \varepsilon$ are positive constants.
	See for instance \cite[\S 5 Lemma 10]{Burgess}. On the other hand,} 

From the approximate functional equation and a dyadic partition of unity, one has
$$L\left(\tfrac 12, \chi \right) \ll_{\varepsilon} M^\varepsilon \sum_N \frac {|S_{\chi} (N) |} {\sqrt N} \left(1 + \frac {N} {\sqrt M}\right)^{-A}, $$
where  $A > 0$,  $N$ ranges over $2^{\nu/2}$ for $\nu = -1, 0 , 1, 2 ...$, and $ S_{\chi} (N) $ is the smooth character sum associated to $\chi$ defined in \eqref{2eq: smooth character sum} for some weight function $W$.
%
%\begin{equation*} %\label{2eq: smooth character sum}
%	S_{\chi} (N) := \sum_{n =1}^\infty \chi (n) W \left(\frac n N \right),
%\end{equation*}
%and $W $ is a smooth function supported on the interval $[1, 2]$ and satisfying $W^{(j)} (x) \ll_j 1$. % We remark that, by means of a dyadic partition of unity and partial summation, \eqref{2eq: character sum} and \eqref{2eq: smooth character sum} are actually equivalent.

The contribution from those $N \gg M^{1/2 + \varepsilon}$ is made negligible by choosing $A$ above to be sufficiently large. Trivially, $|S_{\chi} (N)| \ll N$, and therefore if $N \ll M^{1/2 - \delta}$ for a given $\delta > 0$ then $ {|S_{\chi} (N) |}/ {\sqrt N} \ll M^{1/4 - \delta/2} $. Hence we are left with
\begin{equation}\label{2eq: sum to L}
L\left(\tfrac 12, \chi \right) \ll_{\varepsilon } \left(M^{1/4-\delta/2} + \max_{  M^{1/2 - \delta} \ll N \ll M^{1/2 + \varepsilon}} \frac {|S_{\chi} (N) |} {\sqrt N} \right) M^{\varepsilon}. 
\end{equation}
Thus subconvexity bounds will now follow if one is able to non-trivially bound $S_{\chi} (N)$  for $M^{1/2 - \delta} \ll N \ll M^{1/2 + \varepsilon}$.

\subsection{A modified $\delta$-symbol method}
\label{mod-delta}
One of our main analytic tools for the proof of Theorem \ref{thm: character sum} will be a version of the circle method introduced in \cite{DFI} and \cite{H}. We start with a smooth approximation of the $\delta$-symbol as described in \cite{H}.

\begin{lemma}\label{lemmacm}
	For any $Q > 1$ there is a positive constant $c_Q$, and a smooth function $h(x,y)$ defined on $(0,\infty)\times\mathbb (-\infty, \infty)$, such that
	\begin{align}
	\label{cm}
	\delta(n,0)=\frac{c_Q}{Q^2}\sum_{q=1}^{\infty}\;\sideset{}{^{\star}}\sum_{a \bmod{q}}e\left(\frac{an}{q}\right)h\left(\frac{q}{Q},\frac{n} {Q^2}\right).
	\end{align}
	The constant $c_Q$ satisfies $c_Q=1+O_A(Q^{-A})$ for any $A>0$. Moreover $h(x,y)\ll 1/x $ for all $y$, and $h(x,y)$ is non-zero only for $x\leq\max\{1,2|y|\}$. 
\end{lemma} 
The smooth function $h(x,y)$ satisfies (see \cite{H}) 
\begin{eqnarray}
x^{i} \frac{\partial^i h}{\partial x^i}(x,y)\ll_i \frac 1 x & \textnormal{and} & \frac{\partial h}{\partial y}(x,y)=0 \label{hbound1}
\end{eqnarray}
for $x\leq 1$ and $|y|\leq x/2$. Also for $|y|>x/2$, we have
\begin{equation}
x^i y^j \frac{\partial^{i + j} h}{\partial x^i \partial y^j}  (x,y)\ll_{i,j} \frac 1 x. \label{hbound2}
\end{equation}

Our variant of the  $\delta$-method makes use of the following observation,
$$
\delta(n,0) =\mathbbm{1}_{K|n} \; \delta(n/K,0),
$$
where $K$ is a positive integer and $\mathbbm{1}_{K|n}$ is equal to $1$ or $0$ according as $K|n$ or not.
From this  \eqref{cm} may be written in the following form upon detecting the condition $K|n$ with additive characters.
\begin{align}
\label{cm2}
\delta(n,0)=\frac{c_Q}{KQ^2} \sum_{q=1}^{\infty}\;\sideset{}{^{\star}}\sum_{a \bmod{q}}\sum_{b \bmod{K}}  e\left(\frac{an}{qK}\right) e\left(\frac{bn}{K}\right) h\left(\frac{q}{Q},\frac{n} {KQ^2}\right).
\end{align}

\delete{\begin{rem} Note that for any non-zero integer $K$, one may write
		$$
		\delta(n,0) =\mathbbm{1}_{K|n} \; \delta(n/K,0) .
		$$
		Therefore, \eqref{cm} may be written in the following form upon detecting the condition $K|n$ with additive characters.
		\begin{align}
		\label{cm2}
		\delta(n,0)=\frac{c_Q}{KQ^2} \sum_{q=1}^{\infty}\;\sideset{}{^{\star}}\sum_{a \bmod{q}}\sum_{b \bmod{K}}  e\left(\frac{an}{qK}\right) e\left(\frac{bn}{K}\right) h\left(\frac{q}{Q},\frac{n} {KQ^2}\right).
		\end{align}
	\end{rem}
	We will apply \eqref{cm2} with $K$ prime in the proof of Theorem \ref{thm: character sum}.  One may compare this with Harcos' application of Jutila's version of the circle method in \cite{Harcos:shifted} for studying shifted convolution sums. Equation \eqref{cm2} introduces the structure of Ramanujan sums with moduli divisible by $K$ while maintaining an explicit error term. In other words, breaking apart the sum over $b$ and $q$ in \eqref{cm2} according to divisibility by prime $K$ one has
	\begin{align}
	\label{mod-delta-sym-meth}
	\delta(n,0)=\frac{c_Q}{KQ^2} \sum_{\substack{q=1\\ (q,K)=1}}^{\infty}\;\sideset{}{^{\star}}\sum_{c \bmod{qK}} e\left(\frac{cn}{qK}\right)  h\left(\frac{q}{Q},\frac{n} {KQ^2}\right)+E(n,K)
	\end{align}
	where 
	$$
	E(n,K)= \frac{c_Q}{KQ^2} \sum_{q=1}^{\infty}\;\sideset{}{^{\star}}\sum_{a \bmod{q}}\sum_{\substack{b \bmod{K}\\ q b \equiv 0 \bmod{K}}}  e\left(\frac{an}{qK}\right) e\left(\frac{bn}{K}\right) h\left(\frac{q}{Q},\frac{n} {KQ^2}\right).
	$$
	which trivially satisfies the bound $E(n,K) \ll 1/K $. }

\subsection{Deligne Bound for Character Sums}\label{sec:deligne}

Suppose $p$ is a prime and $\chi$ is a primitive character modulo $p$. Let $m ,  n \in \mathbb F_p$. Consider the sum 
\begin{equation}\label{2eq: character sum}
\mathfrak{S}_{\chi} (m , n) =\sum_{\substack{x_1\in\mathbb{F}_p^\star}}\chi\left( x_1 \right) \bar\chi\left(m + x_1 \right)
e\left(\frac{n x_1}{p}\right).
\end{equation}
Clearly,  we have $\mathfrak{S}_{\chi} (0 , 0) = p-1$. In the case  $m = 0$ and $n \neq 0$ the character sum  reduces to 
\begin{align*}
\mathfrak{S}_{\chi} (0 , n) = \sum_{\substack{x_1\in\mathbb{F}_p^\star }}e\left(\frac{n x_1}{p}\right) = -1.
\end{align*}
In the case $ n = 0 $ and $m \neq 0$ we have
\begin{equation*}
\mathfrak{S}_{\chi} (m , 0) = \sum_{\substack{x_1\in\mathbb{F}_p^\star}} \bar\chi\left(m \overline {x}_1 + 1 \right) = - 1.
\end{equation*}
Finally, we suppose $ m n \neq 0 $. In this case we will use the following relation
$$
\bar{\chi}(a)g_\chi 
=\sum_{b\in\mathbb{F}_p^\star}\chi(b)e\left(\frac{a b}{p}\right) ,
$$
which holds for any $a\in\mathbb{F}_p$. Here 
$$
g_\chi=\sum_{a\in\mathbb{F}_p^\star}\chi(a)e\left(\frac{a}{p}\right)
$$
is the Gauss sum associated with the character $\chi$. Using this relation we rewrite the above character sum as
\begin{equation}
\label{char-sum-1}
\mathfrak{S}_{\chi} (m , n)=\frac{1}{g_\chi  } \sum_{\substack{\mathbf{x}\in \mathbb{F}_p^{\star 2}}}\chi(x_1) {\chi}(x_2)
e\left(\frac{f(\mathbf{x})}{p}\right),
\end{equation}
where
$$
f(\mathbf{x})= n  x_1 +  m x_2 +x_1 x_2.
$$\\

\delete{
	\textcolor{red}{remove below. Unfortunately trivial estimate}
	We first state a simple non-trivial estimate for $\mathfrak{S}_{\chi} (m , n)$ obtained via a second moment method.  Consider the sum
	\begin{equation}
	\label{char-second-moment}
	\frac{1}{\varphi(p)} \sum_{\chi \bmod{p}} \left| \sum_{\substack{\mathbf{x}\in \mathbb{F}_p^{\star 2}}}\chi(x_1) {\chi}(x_2)
	e\left(\frac{f(\mathbf{x})}{p}\right)\right|^2.
	\end{equation}
	Opening the square and summing over the characters $\chi$ modulo $p$ we obtain
	$$
	\mathop{\sum\sum}_{\substack{\mathbf{x, y}\in \mathbb{F}_p^{\star 2} \\ x_1 x_2 \equiv y_1 y_2 \bmod{p}}} e\left(\frac{n (x_1-y_1) + m(x_2-y_2)}{p}\right)
	$$
	which equals
	$$
	\mathop{\sum\sum\sum}_{x_1, \; y_1, \; x_2 \in  \mathbb{F}_p^{\star}} e\left(\frac{n (x_1-y_1) + x_2\; m (1-x_1 \bar{y_1})}{p}\right)
	$$
	and is easily bounded by $p^2$.  Therefore, we obtain the bound
	$$
	\sum_{\substack{\mathbf{x}\in \mathbb{F}_p^{\star 2}}}\chi(x_1) {\chi}(x_2)
	e\left(\frac{f(\mathbf{x})}{p}\right) \ll p^{3/2}
	$$
	so that
	\begin{equation}
	\label{char-second-moment-bound}
	\mathfrak{S}_{\chi} (m , n) \ll p
	\end{equation}
	in the case of $m\neq 0$. \textcolor{red}{Although seemingly non-trivial for \eqref{char-sum-1}, still trivial for \eqref{2eq: character sum}}
	\\
}

In order to obtain a non-trivial estimate, we analyse the sum \eqref{char-sum-1} using Deligne's work as has been developed in \cite{Fu:Deligne}. Let us briefly recall the main result of \cite{Fu:Deligne} concerning sums of the form 
\begin{align*}
\sum_{\substack{\mathbf{x}\in \mathbb{F}_p^{\star r} }}\chi_1(x_1)\chi_2(x_2)\dots \chi_r(x_r)
\psi\left(f(\mathbf{x})\right),
\end{align*}
where $\psi$ is a non-trivial additive character modulo $p$. Here 
$$
f(\mathbf{x})=\sum_{\mathbf{i}\in \mathbb{Z}^r}a_\mathbf{i}\mathbf{x}^\mathbf{i}
$$ is a Laurent polynomial with coefficients $a_\mathbf{i}\in \mathbb{F}_p$. Let $\Delta_\infty(f)$ be the Newton polyhedron associated with $f$. This is given by the convex hull in $\mathbb{R}^r$ of the set
$$
\{\mathbf{i}\in\mathbb{Z}^r:a_\mathbf{i}\neq 0\}\cup\{\mathbf{0}\}.
$$  
The Laurent polynomial $f$ is said to be non-degenerate with respect to $\Delta_\infty(f)$ if for any face $\tau$ of $\Delta_\infty(f)$ not containing the origin, the locus 
\begin{align*}
\frac{\partial f_\tau}{\partial x_1}=\dots=\frac{\partial f_\tau}{\partial x_r}=0
\end{align*}
in the torus $\mathbb{T}_{{\mathbb{F}}_p}^r = {\mathbb{F}}_p^{\star r}$ is empty, where $f_\tau $ denotes the sub-polynomial
$$
f_\tau(\mathbf{x})=\sum_{\mathbf{i}\in \tau}a_\mathbf{i}\mathbf{x}^\mathbf{i}.
$$
If $\dim \Delta_\infty(f)=r$ and $f$ is non-degenerate with respect to $\Delta_\infty(f)$, then we have
\begin{align}
\label{char-sum-del}
\sum_{\substack{\mathbf{x}\in\mathbb{F}_p^{\star r}}}\chi_1(x_1)\chi_2(x_2)\dots \chi_r(x_r)
\psi\left(f(\mathbf{x})\right)\ll p^{r/2},
\end{align}
where the implied constant is independent of $p$.\\

Let us now return to the special case of \eqref{char-sum-1} with $nm \neq 0$.  The Newton polyhedron $\Delta_\infty(f)$ of $f$ is given by the convex hull of 
$$
\{\mathbf{0},\mathbf{e}_1,\mathbf{e}_2, \mathbf{e}_1+\mathbf{e}_2 \},
$$ 
which is $2$ dimensional.
Here $\mathbf{e}_1=(1,0)$ and $\mathbf{e}_2 = (0, 1)$ are the standard basis vectors. 
We have
$$\frac{\partial f}{\partial x_1}  = n  +x_2, \hskip 10 pt \frac{\partial f}{\partial x_2} = m +x_1 .$$
Let $g$ be a sub-polynomial of $f$ such that the equations $\partial g/\partial x_1 (\mathbf {x}) = \partial g/\partial x_2 (\mathbf {x}) = 0$ are solvable on $ {\mathbb{F}}_p^{\star 2}$. 
It is easy to verify that one must have $g = 0$ or $g=f$. It is clear that
neither $0$ nor $f$ is equal to $f_{\tau}$ for any face $\tau$ of $\Delta_{\infty} (f)$ not containing the origin. This proves that $f$ is non-degenerate with respect to $\Delta_\infty(f)$. Using \eqref{char-sum-del} along with the expression \eqref{char-sum-1} of $\mathfrak{S}_{\chi} (m , n)$, we obtain
$$  \mathfrak{S}_{\chi} (m , n) \ll \sqrt p.$$  
We have arrived at the following Lemma.

\begin{lemma}\label{lem: character sum}
	Let $p$ be a prime and $\chi$ be a primitive character modulo $p$. For $m ,  n \in \mathbb F_p$ define the character sum $\mathfrak{S}_{\chi} (m , n)$ by \eqref{2eq: character sum}. Then we have
	\begin{itemize}
		\item[-] $  \mathfrak{S}_{\chi} (0 , 0) = p-1$, 
		\item[-] $ \mathfrak{S}_{\chi} (m , n) = -1$ if $n m = 0$ and either $m \neq 0$ or $n \neq 0$, and
		\item[-] $  \mathfrak{S}_{\chi} (m , n) \ll \sqrt p$ if $m n \neq 0$.
	\end{itemize}
\end{lemma}

\section{Proof of Theorem~\ref{thm: character sum} and Corollary~\ref{cor: subconvexity}}

Let $M_i$, with $i=1,2,3$, be three distinct primes, and set $M=M_1M_2M_3$. Let $\chi_i$ be primitive characters modulo $M_i$ and set $\chi:=\chi_1 \chi_2 \chi_3$. Suppose 
$W $ is a real-valued smooth function on $(0, \infty)$ supported in $[1,2]$ and satisfying 
$$
W^{(j)}(x)\ll_j 1.
$$ 
We shall consider the smooth character sum 
\begin{align*}
%\label{smooth-char-sum}
S_\chi(N)  =\sum_{n\in\mathbb{Z}} \chi_1 \chi_2 \chi_3(n)W\left(\frac{n}{N}\right)
\end{align*}
when $M_1 \ll  N\ll M_1 \min\big\{M_2^{2/3}, M_3^{2}\big\}$ (conditions which arise in the course of the proof). Our goal is to establish a non-trivial bound which will be used in application to the subconvexity problem.

\subsection{Applying the $\delta$-method}

We first write
$$
S_\chi(N)  =\mathop{\sum\sum}_{n, m\in\mathbb{Z}} \chi_1 \chi_2(n)\chi_3(m) \delta(n-m,0)W\left(\frac{n}{N}\right)V\left(\frac{m}{N}\right),
$$
where $V$ is a smooth function with support $[1/2,3]$ and such that $V(x)=1$ for $x\in [1,2]$ with $V^{(j)}(x)\ll_j 1$.

We apply the modified $\delta$-symbol method described in Section~\ref{mod-delta}, with the divisibility modulus $K = M_1$  and $ Q=\sqrt{N/M_1}. $
From \eqref{cm2} we get
\begin{align*}
S_\chi(N)  = &\frac{c_Q}{N}\mathop{\sum\sum}_{n, m\in\mathbb{Z}} \chi_1 \chi_2(n)\chi_3(m) \sum_{q=1}^{\infty} \ \sideset{}{^{\star}}\sum_{a \bmod{q}}\sum_{b \bmod{M_1}}  \\
&
e\left(\frac{a(n-m)}{qM_1}\right)e\left(\frac{b(n-m)}{M_1}\right)  h\left(\frac{q}{Q},\frac{n-m} {N}\right)W\left(\frac{n}{N}\right)V\left(\frac{m}{N}\right).
\end{align*}
In order to have enough points of summation, it is required that 
\begin{equation}
\label{size condition 0}
M_1 \ll N.
\end{equation}

Estimating trivially at this stage we get $S_\chi(N) \ll N^2$. So our job is to save more than $N$.

\subsection{Poisson summation} 

%We will now analyse $S   $. 

\subsubsection{Poisson summation in the $m$-sum.} 

Poisson summation over $m$ gives
\begin{align*}
\mathop{\sum}_{m\in \mathbb{Z}} \chi_3(m) e\left(- \frac{am}{qM_1} - \frac{b m }{M_1}\right)& h\left(\frac{q}{Q},\frac{n-m} {N}\right)V\left(\frac{m}{N}\right)\\
&=\frac{N}{qM_1M_3}\sum_{m\in\mathbb{Z}}\mathfrak{C}(m, q, a, b)\EuScript{I}(m, n, q)
\end{align*}
where the character sum is given by
\begin{align*}
\mathfrak{C}(m, q, a, b) = \sum_{c\bmod{qM_1M_3}}\chi_3(c) e\left(- \frac{ac}{qM_1} - \frac{bc}{M_1} +\frac{mc}{qM_1M_3}\right),
\end{align*}
and the integral is
\begin{align*}
\EuScript{I}(x, v, q)=\int_{\mathbb{R}} h\left(\frac{q}{Q},\frac{v} {N}-u\right)V\left(u\right)e\left(-\frac{Nx u}{qM_1M_3}\right)\mathrm{d}u.
\end{align*}
Applying integration-by-parts and the bounds from \eqref{hbound2}, we see that if $|m|\gg QM_1M_3M^\varepsilon/N$ then the integral is negligibly small (i.e. $O_A(M^{-A})$ for any $A>0$).
%\begin{align*}
%&\mathop{\sum}_{m\in \mathbb{Z}} \chi_3(m)e\left(-\frac{cm}{M_1q}\right)h\left(\frac{q}{Q},\frac{n-m} {N}\right)V\left(\frac{m}{N}\right)\\
%&=\frac{N}{qM_1M_3}\sum_{|m|\ll qM_1M_3/N}\sum_{a\bmod{qM_1M_3}}\chi_3(a)e\left(-\frac{aa}{M_1q}\right)e\left(-\frac{ba}{M_1}\right)e\left(\frac{ma}{qM_1M_3}\right).
%\end{align*}
%Want $N<QM_1M_3$, i.e. $M_2<M_1M_3^3$ (so that we have non-zero frequencies). 

We impose the restriction that $Q\ll M_3$ with a sufficiently small implied constant. This is equivalent to having
\begin{align}
\label{size-condition-1}
N \ll M_1 M_3^2 %M_2\ll M_1M_3^3.
\end{align}
Under this condition we have $(q, M_3) = 1$, and hence the character sum $\mathfrak{C}(m, q, a, b)$ splits as
\begin{align*}
\sum_{c\bmod{qM_1}}e\left(-\frac{ac}{qM_1 } -\frac{bc}{M_1} + \frac{\bar{M_3}mc}{qM_1}\right)\sum_{d\bmod{M_3}}\chi_3(d) e\left(\frac{\overline{qM_1} m d}{M_3}\right).
\end{align*}
It vanishes save for $(a+qb) M_3\equiv m\bmod{qM_1}$ and $(m, M_3) = 1$ in which case we have
$$
qM_1\varepsilon_{3}\sqrt{M_3} \chi_3(qM_1 ) \overline \chi_3 (m).
$$
(Here $\varepsilon_i$ stands for the sign of the Gauss sum associated with the character $\chi_i$, i.e. $g_{\chi_i} = \varepsilon_i \sqrt {M_i}$.) Observe that the congruence condition above implies $(m, q) = 1$.
It follows that
\begin{align*}
S_\chi(N)   =&\frac{\varepsilon_3 c_Q}{\sqrt{M_3}} \sum_{q \ll Q }\ \sum_{\substack{ |m|\ll QM_1M_3M^\varepsilon/N \\ (m,qM_3)=1}} \chi_3(qM_1 ) \overline \chi_3 (m) \mathop{\sum}_{n\in\mathbb{Z}} \chi_1 \chi_2(n) \\
&\mathop{\sideset{}{^\star}\sum_{a\bmod{q}}\ \sum_{b\bmod{M_1}}}_{(a+qb)M_3\equiv m\bmod{qM_1}}  e\left(\frac{an}{qM_1} + \frac {bn} {M_1} \right)\EuScript{I}(m, n, q) W\left(\frac{n}{N}\right)+O_A(M^{-A}).
\end{align*}
We thus need to consider the sum
\begin{equation*}
\begin{split}
\frac{1}{\sqrt{M_3}}\sum_{q\ll Q } \ &\sum_{\substack{ |m|\ll QM_1M_3M^\varepsilon/N \\ (m,qM_3)=1}} \chi_3(q ) \overline \chi_3 (m) \\
&\mathop{\sum}_{n\in\mathbb{Z}} \chi_1 \chi_2(n) e\left(\frac{\bar{M}_3mn}{qM_1}\right)\EuScript{I}(m, n, q)W\left(\frac{n}{N}\right).
\end{split}
\end{equation*}
At this stage trivial estimation gives
$$
S_\chi(N)   \ll N\sqrt{M_3}M^\varepsilon,
$$
and it remains to save more than $\sqrt{M_3}$.\\

\subsubsection{Poisson summation in the $n$-sum} 

Next, we apply the Poisson summation formula on the sum over $n$. This gives
\begin{align*}
\mathop{\sum}_{n\in\mathbb{Z}} \chi_1\chi_2(n)e\left(\frac{\bar{M}_3mn}{qM_1}\right)& \EuScript{I}(m, n, q) W\left(\frac{n}{N}\right)\\
&=\frac{N}{q M_1M_2 }\sum_{n\in\mathbb{Z}} \mathfrak{C} (m, n, q) \EuScript{J}(m, n, q),
\end{align*}
where the character sum is given by
\begin{align*}
\mathfrak{C} (m, n, q) = \sum_{a\bmod{M_1M_2q}}\chi_1\chi_2(a)e\left(\frac{\bar{M_3}ma}{q M_1 } + \frac{na}{q M_1M_2}\right)
\end{align*}
and the integral is
\begin{align*}
\begin{split}
\EuScript{J}(x, y, q) & = \int_{\mathbb{R} } \EuScript{I}(x, Nv, q)  W(v) e\left(-\frac{ Nvy}{qM_1M_2}\right) \mathrm{d}v\\
& = \iint_{\mathbb{R}^2} h\left(\frac{q}{Q},v-u\right)W(v)V\left(u\right) e\left(-\frac{ Nux}{qM_1M_3}-\frac{ Nv y}{qM_1M_2}\right)\mathrm{d}u\mathrm{d}v.
\end{split}
\end{align*}
By repeated integration-by-parts we get that the tail $|n|\gg QM_1M_2M^\varepsilon/N$ makes a negligible contribution to the sum.
%*****Edited so far**** (2Aug14.RM)*****

We impose the restriction that $Q\ll M_2$ with a sufficiently small implied constant. This is equivalent to having
\begin{align}
\label{size-condition-2}
N \ll M_1 M_2^2 %M_3\ll M_1M_2^3.
\end{align}
Under this condition  $(q, M_2) = 1$, and therefore the character sum splits as
\begin{align*}
\sum_{a\bmod{qM_1}}\chi_1(a)e\left(\frac{\bar{M_3}ma}{qM_1}\right)e\left(\frac{\bar{M_2}na}{qM_1}\right) \sum_{b\bmod{M_2}}\chi_2(b)e\left(\frac{\overline{qM_1}nb}{M_2}\right),
\end{align*}
which is
\begin{align*}
\varepsilon_2\sqrt{M_2}\chi_2(qM_1 ) \overline \chi_2 (n) \sum_{a\bmod{qM_1}}\chi_1(a)e\left(\frac{(\bar{M_3}m + \bar{M_2}n)a}{qM_1} \right).
\end{align*}
Suppose $q=q' M_1^r$ with $M_1\nmid q'$. Then the remaining character sum splits into the product
\begin{align*}
\sum_{a\bmod{q'}} & e\left(\frac{\big(\overline{M_3}m + \overline{M_2 }n\big) \overline {M_1^{r+1}} a}{q'}  \right)\\
&\times \sum_{b\bmod{M_1^{r+1}}} \chi_1(b) e\left(\frac{(\overline{M_3}m + \overline{M_2}n) \overline{q'} b}{M_1^{r+1}} \right).
\end{align*}
This product vanishes unless $q'M_1^r|M_2m+M_3n$ in which case we get
\begin{align*}
q' \varepsilon_1\sqrt{M_1} M_1^r\chi_1(q'M_2M_3)\bar\chi_1\left((M_2m+M_3n)/M_1^r\right).
\end{align*}
We conclude that
\begin{align*}
&\mathop{\sum}_{n\in\mathbb{Z}} \chi_1\chi_2(n)e\left(\frac{\bar{M}_3mn}{qM_1}\right)W\left(\frac{n}{N}\right)\mathfrak{I}(m, n, q)=\frac{\varepsilon_1\varepsilon_2 N}{\sqrt{M_1M_2}}\\
& \sum_{\substack{|n|\ll QM_1M_2M^\varepsilon/N\\q|M_2m+M_3n}}\chi_1(q'M_2M_3)\bar\chi_1\left((M_2m+M_3n)/M_1^r\right)\chi_2(qM_1 ) \overline \chi_2 (n)\EuScript{J}(m, n, q) \\
& \hskip 295 pt +O_A(M^{-A}),
\end{align*}
where $r=v_{M_1}(q)$ is the $M_1$-adic valuation of $q$ and $q' = q/M_1^r$. Since $(m, q ) = 1$, we also have $(n, q ) = 1$. Consequently
\begin{equation*}
\begin{split}
& S_\chi(N)   = \frac {\eta c_Q N} {\sqrt M}  \sum_{\substack{q\ll Q }} \chi_1(q') \\
& \mathop{\sum\sum}_{\substack{|n|\ll QM_1M_2M^\varepsilon/N\\|m|\ll QM_1M_3M^\varepsilon/N\\(m,qM_3) =  (n, qM_2) =1\\q|M_2m+M_3n}} \bar\chi_1\left((M_2m+M_3n)/M_1^r\right) \chi_2(q ) \overline \chi_2 (n) \chi_3(q )\overline \chi_3 (m)\EuScript{J}(m, n, q)\\
&\hskip 285 pt +O_A(M^{-A}),
\end{split}
\end{equation*}
with $|\eta|=1$. At this stage, trivial estimation gives
$$
S_\chi(N)   \ll M^{1/2 + \varepsilon},
$$
which is just at the threshold. Any additional saving will yield a non-trivial bound for the character sum.\\

Observe that for $r\geq 1$ we are saving an extra $M_1$ by trivial estimation. Therefore, we just need to focus on the generic case $r=0$. We consider
\begin{equation*}
\begin{split} 
S_0   := & \sum_{\substack{q\ll Q\\ (q,M_1)=1}} \chi_1\chi_2\chi_3(q)\\
& \mathop{\sum\sum}_{\substack{|n|\ll QM_1M_2M^\varepsilon/N\\|m|\ll QM_1M_3M^\varepsilon/N\\(m,qM_3) = (n, q M_2) =1\\q|M_2m+M_3n}} \bar\chi_1\left(M_2m+M_3n\right)\bar\chi_2(n)\bar\chi_3(m)\EuScript{J}(m, n, q).
\end{split}
\end{equation*}
Then
\begin{equation}\label{3eq: r neq 0 error}
S_\chi(N)   = \frac {\eta c_Q N} {\sqrt M} S_0   + O \left(\frac {N M^\varepsilon}  {M_1}\right).
\end{equation}

\subsection{Treatment of $S_0$}

\subsubsection{Applying Cauchy's inequality}
From Cauchy's inequality we get
\begin{equation}\label{3eq: S0 and T}
S_0  \ll M^\varepsilon\sqrt{\frac{QM_1M_2}{N}}\sqrt{T  },
\end{equation}
where $T$ is given by  \small
\begin{equation*}
\sum_{|n|\ll QM_1M_2M^\varepsilon/N}\left|\sum_{\substack {q\ll Q \\ (q, M_1 n)  = 1 }} \chi_1\chi_2\chi_3(q) \mathop{\sum}_{\substack{|m|\ll QM_1M_3M^\varepsilon/N\\(m,qM_3)=1\\q|M_2m+M_3n}}  \bar\chi_1\left(M_2m+M_3n\right)\bar\chi_3(m)\EuScript{J}(m, n, q)\right|^2.
\end{equation*}\normalsize
Here we assumed that $M_3<M_2$, otherwise we would have pulled out the $m$-sum rather than the $n$-sum. Any non-trivial bound for $T$ will yield a non-trivial bound for the character sum. Introducing a dyadic partition of unity for the $n$-sum and opening the absolute square it suffices to consider the following sum
\begin{equation}\label{3eq: open the square}
\mathop{\sum\sum}_{\substack{q_1, q_2\ll Q \\ (q_i, M_1) = 1}}\chi_1\chi_2\chi_3(q_1\bar{q_2}) 
\mathop{\sum\sum}_{\substack{|m_i|\ll QM_1M_3M^\varepsilon/N\\(m_i,q_iM_3)=1}} \chi_3(\overline{m_1} m_2 )\:{T} (m_1, m_2, q_1, q_2 )
\end{equation}
where ${T} (m_1, m_2, q_1, q_2 )$ is given by
\begin{equation*}
\begin{split}
\sum_{\substack{n\in\mathbb{Z}\\n\equiv -M_2\bar{M_3}m_i\bmod{q_i}}} & \bar\chi_1\left(M_2m_1+M_3n\right)\chi_1\left(M_2m_2+M_3n\right) \\
&\times \EuScript{K}(m_1, m_2, n, q_1, q_2 )  U\left(\frac{n}{R}\right),
\end{split}
\end{equation*}
with
$$
\EuScript{K}(x_1, x_2, y, q_1, q_2 )= \EuScript{J} (x_1, y, q_1 ) \overline{\EuScript{J} (x_2, y, q_2)}.
$$
Here $U$ is a suitable smooth function with compact support and $R \ll QM_1M_2M^\varepsilon/N$.  

\subsubsection{The third application of Poisson summation} 
We seek to get cancellation in ${T} (m_1, m_2, q_1, q_2 )$. For this we at least need that the sum has enough points of summation or
$
QM_1M_2  /N\gg Q^2 
$, which is equivalent to
\begin{align}
\label{size-condition-3}
N\ll M_1M_2^{2/3}.
\end{align}  
We now apply Poisson summation to ${T} (m_1, m_2, q_1, q_2 )$ with modulus $q_1q_2M_1$. This gives
$$
{T} (m_1, m_2, q_1, q_2 )=\frac{R}{q_1q_2M_1}\sum_{n\in\mathbb{Z}}\mathfrak{C}(m_1, m_2, n, q_1, q_2) \EuScript{L} ( m_1, m_2, n, q_1, q_2)
$$
where the character sum $\mathfrak{C}(m_1, m_2, n, q_1, q_2)$ is given by
$$
\sum_{\substack{a\bmod{q_1q_2M_1}\\a\equiv -M_2\bar{M_3}m_i \bmod{q_i}}} \bar\chi_1\left(M_2m_1+M_3a\right)\chi_1\left(M_2m_2+M_3a\right)e\left(\frac{na}{q_1q_2M_1}\right),
$$
and the integral is given by
\begin{equation*}
\EuScript{L} ( x_1, x_2, z, q_1, q_2) = \int_{\mathbb R} \EuScript{K}(x_1, x_2, R y, q_1, q_2 ) U(y) e \left( - \frac {R y z} {q_1 q_2 M_1} \right)  \mathrm{d} y.
\end{equation*}
By repeated integration-by-parts we see that the integral $\EuScript{L} ( m_1, m_2, n, q_1, q_2)$ is negligibly small if $|n|\gg Q^2M_1M^\varepsilon/R = N M^\varepsilon/R   $. %or $|n|\gg  \max\{q_1,q_2\} NM^\varepsilon/M_2$. 
Hence
\begin{align*}
& {T} (m_1, m_2, q_1, q_2 ) =  \frac{R}{q_1q_2M_1} \\
&  \sum_{|n|\ll N M^\varepsilon/R}\mathfrak{C}(m_1, m_2, n, q_1, q_2) \EuScript{L} ( m_1, m_2, n, q_1, q_2) +O_A(M^{-A}).
\end{align*}
Using the trivial bound $|\EuScript{L} ( m_1, m_2, n, q_1, q_2) |\ll Q^2/q_1q_2$ which follows from \eqref{hbound1}, we conclude that
\begin{align*}
{T} (m_1, m_2, q_1, q_2 )  
\ll \frac{RQ^2}{(q_1q_2)^2M_1}\sum_{|n|\ll N M^\varepsilon/R}|\mathfrak{C}(m_1, m_2, n, q_1, q_2)|+ M^{-A}.
\end{align*}

\subsubsection{Bounds for $\mathfrak{C}(m_1, m_2, n, q_1, q_2)$ and ${T} (m_1, m_2, q_1, q_2 )$}
Since $ M_1 \nmid q_1q_2 $, the character sum  $\mathfrak{C}(m_1, m_2, n, q_1, q_2)$ splits as
\begin{align*}
& \sum_{\substack{a\bmod{M_1}}} \bar\chi_1\left(M_2m_1+ M_3a\right) \chi_1\left(M_2m_2+ M_3a\right)
e\left(\frac{ \overline{q_1q_2} n a}{M_1}\right)\\
&\times \sum_{\substack{b\bmod{q_1q_2}\\b\equiv -M_2\bar{M_3}m_i\bmod{q_i}}}e\left(\frac{\bar{M_1}n b}{q_1q_2}\right).
\end{align*}
The second sum %vanishes unless $ (q_1, q_2)  |m_2-m_1$ in which case we have 
has bound $(q_1, q_2)$, since there are at most $(q_1, q_2)$ many terms due to the Chinese remainder theorem. %vanishes unless  $(q_1, q_2) | n$ in which case we get a number of norm $(q_1, q_2)$.  
On making the change of indices $x_1 = M_2 m_2 + M_3 a$, one sees that the first sum is equal to%in the non-trivial case $(q_1, q_2) = 1$ 
$$%\mathfrak{C}(m_1, m_2, n, q_1, q_2)= 
\eta\, \mathfrak{S}_{\chi_1} \left(M_2 (m_1 - m_2) \bmod{M_1}, \overline {q_1 q_2 M_3} n  \bmod{M_1}\right)$$ for some $\eta$ with $|\eta| = 1$. Here we recall that the character sum $\mathfrak{S}_{\chi} (m , n)$ is defined by \eqref{2eq: character sum} in  Section~\ref{sec:deligne}.

At this point we need to apply bounds from Section~\ref{sec:deligne}.  In view of Lemma \ref{lem: character sum}, we have the following uniform bound 
\begin{equation*}
\mathfrak{C}(m_1, m_2, n, q_1, q_2) \ll (q_1, q_2) \sqrt { M_1 (n,M_1) }.
\end{equation*}
However, for the zero frequency $n=0$ we shall use the bound
\begin{equation*}
|\mathfrak{C}(m_1, m_2, 0, q_1, q_2)| \leq (q_1, q_2) (m_1 - m_2, M_1).
\end{equation*}
\delete{
	First consider the diagonal contribution $m_1\equiv m_2\bmod{M_1}$. In this case the character sum modulo $M_1$ reduces to 
	\begin{align*}
	\sum_{\substack{a\bmod{M_1}\\a\equiv - M_2\overline{M}_3m_1\bmod{M_1}}} e\left(\frac{\overline{q_1q_2} n a}{M_1}\right),
	\end{align*}
	which is dominated by $\sqrt{M_1}\sqrt{(n,M_1)}$. In the case where $M_1|n$, this bound coincides with the trivial bound $O(M_1)$. However for the zero frequency $n=0$ we need a non-trivial bound. Indeed in this case the character sum modulo $M_1$ reduces to
	\begin{align*}
	&\sum_{\substack{a\bmod{M_1}}}\bar\chi_1\left(m_1M_2+aM_3\right)\chi_1\left(m_2M_2+aM_3\right).
	\end{align*}
	We set $x=m_1M_2+aM_3$ and $y=M_2(m_2-m_1)$. The character sum is then given by
	\begin{align*}
	&\sum_{\substack{x\bmod{M_1}}}\bar\chi_1\left(x\right)\chi_1\left(x+y\right)=
	\sideset{}{^\star}\sum_{\substack{x\bmod{M_1}}}\chi_1\left(1+\bar{x}y\right)\ll (m_1-m_2, M_1).
	\end{align*}
	Next we consider the case where $n$ and $m_1-m_2$ are non-zero modulo $M_1$. In this case Lemma?? yields the bound $O(\sqrt{M_1})$ for the character sum.\\
}
We conclude that ${T} (m_1, m_2, q_1, q_2 )$ is dominated by
$$
\frac{(q_1, q_2) RQ^2}{(q_1q_2)^2M_1}\left( (m_1-m_2, M_1)+\sqrt{M_1}\sum_{0 \neq |n|\ll NM^\varepsilon/R}\sqrt{(n,M_1)}\right).
$$
This gives
\begin{align}
\label{t'-bd}
{T} (m_1, m_2, q_1, q_2 )\ll \frac{(q_1, q_2) RQ^2 M^\varepsilon}{(q_1q_2)^2M_1}  (m_1-m_2, M_1)+ \frac{ (q_1, q_2) Q^2 N M^\varepsilon }{(q_1q_2)^2 \sqrt{M_1} } .
\end{align}

\subsubsection{Bound for $S_0$}
To bound $T$, we will now substitute the bound \eqref{t'-bd} in to \eqref{3eq: open the square} and estimate the remaining sums trivially. To estimate the contribution of the first term in \eqref{t'-bd}, we observe that
\begin{align*}
\frac{RQ^2}{M_1}\mathop{\sum\sum}_{\substack{q_1, q_2\ll Q}}\mathop{\sum\sum}_{\substack{|m_i|\ll QM_1M_3M^\varepsilon/N\\(m_i,q_iM_3)=1   }}\frac{(q_1, q_2)}{(q_1q_2)^2}
%\ll \frac{RQ^2(M_1M_3)^2M^\varepsilon}{M_1N^2}%(DETAILS of CALCULATION)
\ll \frac{\sqrt{M_1}M_2M_3^2}{N^{3/2}}M^\varepsilon 
\end{align*}
and also 
\begin{align*}
RQ^2\mathop{\sum\sum}_{\substack{q_1, q_2\ll Q}}\mathop{\sum\sum}_{\substack{|m_i|\ll QM_1M_3M^\varepsilon/N\\(m_i,q_iM_3)=1\\  M_1|m_2-m_1}}\frac{(q_1, q_2)}{(q_1q_2)^2}
%(HIDE DETAILS)&\ll RQ^2M^\varepsilon\mathop{\sum\sum}_{\substack{q_1, q_2\ll Q}}\frac{M_1M_3\min\{q_1,q_2\}}{N(q_1q_2)^2}\mathop{\sum}_{\substack{|\ell|\ll M_1M_3\max\{q_1,q_2\}M^\varepsilon/N\\M_1|\ell}}\:1\\
%(HIDE DETAILS)&\ll RQ^2M^\varepsilon\mathop{\sum\sum}_{\substack{q_1, q_2\ll Q}}\frac{M_1M_3\min\{q_1,q_2\}}{N(q_1q_2)^2}\left\{1+\frac{M_3\max\{q_1,q_2\}}{N}\right\}\\
%(HIDE DETAILS)&\ll RQ^2\left\{\frac{M_1M_3}{N}+\frac{(M_1M_3)^2}{M_1N^2}\right\}M^\varepsilon\\
\ll \frac{ {M_1}M_2M_3}{ {N}}\left( 1+\frac{M_3}{{\sqrt {N M_1}}}\right)M^\varepsilon.
\end{align*}
To estimate the contribution of the second term in \eqref{t'-bd} towards $T$ we evaluate
\begin{align*}
\frac{Q^2 N} {\sqrt {M_1}} \mathop{\sum\sum}_{\substack{q_1, q_2\ll Q}}\mathop{\sum\sum}_{\substack{|m_i|\ll QM_1M_3M^\varepsilon/N\\(m_i,q_iM_3)=1 \\   }}&\frac{(q_1, q_2)}{(q_1q_2)^2}\ll \sqrt{M_1} M_3^2M^\varepsilon.
\end{align*}
Inserting these bounds in \eqref{3eq: S0 and T} we obtain
\begin{align}\label{3eq: bound of S0}
S_0  \ll \left(\frac{\sqrt{M_1}M_2M_3}{N} + \frac{M_1^{3/4} M_2 \sqrt {M_3}} {N^{3/4}} +\frac{\sqrt{M_1M_2}M_3}{N^{1/4}}\right)M^\varepsilon.
\end{align}

\subsection{Conclusion}
Observe that the first term in \eqref{3eq: bound of S0} absorbs the error term in \eqref{3eq: r neq 0 error}, and therefore we obtain from \eqref{3eq: r neq 0 error} the following bound for $S_\chi(N)$,
\begin{align*}
S_\chi(N)   \ll_\varepsilon \left(M_2^{1/2} M_3^{1/2}+ {M_1^{1/4} M_2^{1/2}  } {N^{1/4}} + M_3^{1/2}{N^{3/4}}\right)M^\varepsilon
\end{align*}
which matches with \eqref{SchiNbound} in Theorem~\ref{thm: character sum}.  Also note that our assumptions \eqref{size condition 0}, \eqref{size-condition-1}, \eqref{size-condition-2} and \eqref{size-condition-3} produced the condition
$$
M_1 \ll N \ll M_1\min \big\{M_2^{2/3}, M_3^{2}\big\}.
$$

Dividing the above bound for $S_\chi(N)$ by $\sqrt{N}$ and returning to \eqref{2eq: sum to L}, we see that $L\left(\tfrac 12, \chi \right)$ is bounded by 
$$
M^{1/4-\delta/2+\varepsilon} + \max_{  M^{1/2 - \delta} \ll N \ll M^{1/2}} \left( \left(\frac{M_2 M_3}{N}\right)^{ {1}/{2}} + \left(\frac{M_1 M_2^{2}} {N}\right)^{ {1}/{4}} + M_3^{ {1}/{2}} N^{ {1}/{4}}\right) M^{\varepsilon}
$$
for any $\varepsilon>0$ provided that $M_1 \ll N \ll M_1\min\big\{M_2^{2/3}, M_3^{2}\big\}$ is satisfied for all $M^{1/2 - \delta} \ll N \ll M^{1/2 + \varepsilon}$.  Therefore, in order to establish Corollary~\ref{cor: subconvexity}, one needs
\begin{equation} \label{maxbound}
\max_{  M^{1/2 - \delta} \ll N \ll M^{1/2}} \max \left\{ \frac{M_2 M_3}{N}, \frac{M_1^{ {1}/{2}} M_2} {N^{ {1}/{2}}}, M_3 N^{ {1}/{2}}, M_1 \right\} \ll M^{1/2-\delta},
\end{equation}
and
$
M^{1/2+\varepsilon} \ll M_1\min\big\{M_2^{2/3}, M_3^{2}\big\}.
$
One can easily verify that the last bound is always satisfied when \eqref{maxbound} is satisfied.  Choosing $N=M^{1/2 - \delta}$ for the first two terms in \eqref{maxbound} and $N=M^{1/2}$ for the third term in \eqref{maxbound}, we obtain the conditions given by the inequalities in \eqref{eq: relations between theta, 3} in Corollary~\ref{cor: subconvexity}.

%    Text of article.

%    Bibliographies can be prepared with BibTeX using amsplain,
%    amsalpha, or (for "historical" overviews) natbib style.
\bibliographystyle{amsalpha}
%    Insert the bibliography data here.
\bibliography{references}

\end{document}